\documentclass[12pt,leqno,draft]{article}

\begin{document}

\title{A brief introduction to Gromov's notion of hyperbolic groups}

\author{Stephen Semmes \\ [10pt]
\textit{dedicated to Leon Ehrenpreis and Mitchel Taibleson}}

\date{}

\maketitle

	The basic reference is \cite{misha2}, and see also \cite{CDP, CP,
GH, misha1, misha3, harpe, ohshika}, for instance.

\subsubsection*{Basic concepts}

	Let $\Gamma$ be a group, and let $F$ be a finite set of
elements of $\Gamma$.  By a \emph{word} over $F$ we mean a formal
product of elements of $F$ and their inverses.  Every word over $F$
determines an element of the group $\Gamma$, simply using the group
operations.  The ``empty word'' is considered a word over $F$, which
corresponds to the identity element of $\Gamma$.

	If $z$ is a word over $F$, then the length of $z$ is denoted
$L(z)$ and is the number of elements of $F$ such that they or their
inverses are used in $z$, counting multiplicities.  A word $z$ is said
to be irreducible if it does not contain an $\alpha \in F$ next to
$\alpha^{-1}$, i.e., so that all obvious cancellations have been made.
If a word $z$ over $F$ corresponds to the identity element of $\Gamma$,
then $z$ is said to be trivial.

	A finite subset $F$ of a group $\Gamma$ is a set of
\emph{generators} of $\Gamma$ if every element of $\Gamma$ corresponds
to a word over $F$.  A group is said to be \emph{finitely-generated}
if it has a finite set of generators.  Let us make the convention that
a generating set $F$ of a group $\Gamma$ should not contain the
identity element of $\Gamma$.

	Suppose that $\Gamma$ is a group and that $F$ is a finite set
of generators of $\Gamma$.  The \emph{Cayley graph} associated to
$\Gamma$ and $F$ is the graph consisting of the elements of $\Gamma$
as vertices with the provision that $\gamma_1$, $\gamma_2$ in $\Gamma$
are adjacent if $\gamma_2 = \gamma_1 \, \alpha$, where $\alpha$ is an
element of $F$ or its inverse.  Thus this relation is symmetric in
$\gamma_1$ and $\gamma_2$.

	A finite sequence $\theta_0, \theta_2, \ldots, \theta_k$ of
elements of $\Gamma$ is said to define a \emph{path} if $\theta_j$,
$\theta_{j+1}$ are adjacent in the Cayley graph for each $j$, $0 \le j
\le k-1$.  The \emph{length} of this path is defined to be $k$.  We
include the degenerate case where $k = 0$, so that a single element of
$\Gamma$ is viewed as a path of length $0$.

	If $\phi$, $\psi$ are elements of $\Gamma$, then the
\emph{distance} between $\phi$ and $\psi$ is defined to be the
shortest length of a path that connects $\phi$ to $\psi$.  In
particular, note that for any two elements $\phi$, $\psi$ in $\Gamma$
there is a path which starts at $\phi$ and ends at $\psi$.  To see
this, one can write $\psi$ as $\phi \, \beta$ for some $\beta$ in
$\Gamma$, and then use the assumption that $\Gamma$ is generated by
$F$ to obtain a path from $\phi$ to $\psi$ one step at a time.

	These definitions are invariant under \emph{left} translations
in $\Gamma$.  In other words, if $\delta$ is any fixed element of
$\Gamma$, then the tranformation $\gamma \mapsto \delta \, \gamma$ on
$\Gamma$ defines an automorphism of the Cayley graph, and it also
preserves distances between elements of $\Gamma$.  This follows
from the definitions, since the Cayley graph was defined in terms of
right-multiplication by generators and their inverses.

	A basic fact is that this definition of distance does not
depend too strongly on the choice of generating set $F$, in the sense
that if one has another finite generating set, then the two distance
functions associated to these generating sets are each bounded by a
constant multiple of the other.  This is not difficult to check, by
expressing each generator in one set as a finite word over the other
set of generators.  There are only a finite number of these expressions,
so that their maximal length is a finite number.

	Let us continue with the assumption that we have a fixed
generating set $F$ for the group $\Gamma$.  Suppose that $R$ is a
finite set of words over $F$.  We say that $R$ is a set of
\emph{relations} for $\Gamma$ if every element of $R$ is a trivial
word.  The inverses of elements of $R$ are also then trivial words, as
well as conjugates of elements of $R$.  That is, if $r$ is an element
of $R$ and $u$ is any word over $F$, then $u \, r \, u^{-1}$ is the
conjugate of $r$ by $u$, and it is a trivial word since $r$ is.
Products of conjugates of elements of $R$ and their inverses are
trivial words too, as well as words obtained from these through
cancellations, i.e., by cancelling $\alpha \, \alpha^{-1}$ and
$\alpha^{-1} \, \alpha$ whenever $\alpha$ is an element of $F$.  The
combination of $F$ and a set $R$ of relations defines a
\emph{presentation} of $\Gamma$ if every word over $F$ which
corresponds to the identity element of $\Gamma$ can be obtained in
this manner.  The empty word is viewed as being equal to the empty
product of relations, so that it is automatically included.  A group
$\Gamma$ is said to be \emph{finitely-presented} if there is a
presentation with a finite set of generators and a finite set of
relations.  For instance, if $\Gamma$ is the free group with
generators in $F$, then one can take $R$ to be the set consisting of
the empty word, and this defines a presentation for $\Gamma$.

	Let us call a word over $F$ \emph{trivial} if it corresponds
to the identity element of $\Gamma$.  Suppose that $w$ is a trivial
word, with
\begin{equation}
	w = \beta_1 \beta_2 \cdots \beta_n,
\end{equation}
where each $\beta_i$ is an element of $F$ or an inverse of an element
of $F$.  This leads to a path $\theta_0, \theta_1, \ldots, \theta_n$,
where $\theta_0$ is the identity element of $\Gamma$ and $\theta_j$ is
equal to $\beta_1 \beta_2 \cdots \beta_j$ when $j \ge 1$.  Because $w$
is a trivial word, $\theta_n$ is also the identity element in
$\Gamma$, which is to say that this path is a loop that begins and
ends at the identity element.

	Fix a finite set $R$ of relations, so that $F$ and $R$ give a
presentation for $\Gamma$.  Let $w$ be a trivial word over $F$ which
is also irreducible.  Define $A(w)$ to be the smallest nonnegative
integer $A$ for which there exist relations $r_1, r_2, \ldots, r_k$ in
$R$, integers $b_1, b_2, \ldots, b_k$, and words $u_1, u_2, \ldots,
u_k$ over $F$ such that the expression
\begin{equation}
  u_1 r_1^{b_1} u_1^{-1} u_2 r_2^{b_2} u_2^{-1} \cdots u_k r_k^{b_k} u_k^{-1}
\end{equation}
can be reduced to $w$ after cancellations as before,
\begin{equation}
	\sum_{j=1}^k L(u_j) \le A,
\end{equation}
and
\begin{equation}
	\sum_{j=1}^k |b_j| \, L(r_j)^2 \le A.
\end{equation}
Here if $z$ is a word over $F$ and $b$ is an integer, then $z^b$ is
defined in the obvious manner, by simply repeating $z$ $b$ times when
$b \ge 0$, or repeating $z^{-1}$ $-b$ times when $b < 0$.  A
representation of this type for $w$ necessarily exists, since $F$ and
$R$ give a presentation for $\Gamma$.

	The group $\Gamma$ is said to be \emph{hyperbolic} if there
is a nonnegative real number $C_0 \ge 0$ so that
\begin{equation}
	A(w) \le C_0 \, L(w)
\end{equation}
for all irreducible trivial words $w$.  The property of hyperbolicity
does not depend on the choice of finite presentation for $\Gamma$, and
in fact there are other definitions for which one only needs to assume
that $\Gamma$ is finitely generated, and the existence of a finite
presentation is then a consequence.  This characterization of
hyperbolicity is discussed in Section 2.3 of \cite{misha2}.  Some
examples of hyperbolic groups are finitely-generated free groups and
the fundamental groups of compact connected Riemannian manifolds
without boundary and strictly negative curvature.  In particular,
this includes the fundamental group of a closed Riemann surface
with genus at least $2$.

\subsubsection*{Spaces of homogeneous type}

	Let us digress now a bit and review some notions from
real-variable harmonic analysis.  Let $M$ be a nonempty set.
A nonnegative real-valued function $d(x,y)$ on the Cartesian
product $M \times M$ is said to be a \emph{quasimetric} if 
$d(x,y) = 0$ exactly when $x = y$, $d(x,y) = d(y,x)$ for all
$x, y \in M$, and 
\begin{equation}
	d(x,z) \le C \Bigl(d(x,y) + d(y,z)\Bigr)
\end{equation}
for some positive real number $C$ and all $x, y, z \in M$.
If this last condition holds with $C = 1$, then $d(x,y)$ is
said to be a \emph{metric} on $M$.

	If $d(x,y)$ is a quasimetric on $M$ and $a$ is a positive real
number, then $d(x,y)^a$ is also a quasimetric on $M$.  If $d(x,y)$ is
a metric on $M$ and $a$ is a positive real number such that $a \le 1$,
then $d(x,y)^a$ is a metric on $M$ too.  These statements are not
difficult to verify.  There is a very nice result going in the other
direction, which states that if $d(x,y)$ is a quasimetric on $M$,
then there are positive real numbers $C'$, $\delta$ and a metric
$\rho(x,y)$ on $M$ such that
\begin{equation}
	C'^{-1} \, \rho(x,y)^\delta \le d(x,y) \le C' \, \rho(x,y)^\delta
\end{equation}
for all $x, y \in M$.  See \cite{MS1}.

	If $d(x,y)$ is a quasimetric on $M$ and $f$ is a real-valued
function on $M$, then $f$ is said to be \emph{Lipschitz} if there
is a nonnegative real number $L$ such that
\begin{equation}
	|f(x) - f(y)| \le L \, d(x,y)
\end{equation}
for all $x, y \in M$.  In general, for a quasimetric, there may not be
any nonconstant Lipschitz functions.  This is the case when $M = {\bf
R}^n$ equipped with the quasimetric $d(x,y)^a$ with $a > 1$, for which
any Lipschitz function would have to have first derivatives equal to
$0$ everywhere.  However, if $d(x,y)$ is a metric, then $f_p(x) =
d(x,p)$ satisfies the Lipschitz condition with $L = 1$ for all $p$ in
$M$.  This can be checked using the triangle inequality.

	If $d(x,y)$ is a quasimetric on $M$, then one has many of the
same basic notions as for a metric, such as convergence of sequences,
open and closed sets, dense subsets, and so on.  One should be a bit
careful with some of the standard results, since for instance it is
not so clear that an open ball defined using a quasimetric is an open
set, as in the situation of ordinary metrics.  At any rate, it still
makes sense to say that $M$ is separable with respect to a quasimetric
if it has a subset which is at most countable and also dense, and one
can define the topological dimension for $M$ as in \cite{HW}.  The
diameter of a subset can be defined in the usual manner using the
quasimetric, and this permits one to define the Hausdorff dimension of
a nonempty subset of $M$.  A famous result about metric spaces is that
the topological dimension is always less than or equal to the
Hausdorff dimension.  See Chapter VII of \cite{HW}.  This does not
work for quasimetrics in general, and it cannot possibly work.  For if
$(M, d(x,y))$ is a quasimetric space with Hausdorff dimension $s$ and
$a$ is a positive real number, then $(M, d(x,y)^a)$ has Hausdorff
dimension $s / a$, while the topological dimension of $(M, d(x,y)^a)$
is the same as that of $(M, d(x,y))$.

	A quasimetric space $(M, d(x,y))$ is said to have the
\emph{doubling property} if there is a positive real number $C_1$ so
that every open ball $B(x,r) = \{y \in M : d(x,y) < r\}$ in $M$ of
radius $r$ can be covered by a family of at most $C_1$ open balls of
radius $r/2$.  By iterating this condition one obtains that for each
positive integer $l$ and each open ball $B(x,r)$ there is a family of
at most $C_1^l$ open balls of radius $2^{-l} r$ which covers $B(x,r)$.
This is a kind of condition of polynomial growth; if one chooses
$\alpha \ge 0$ so that $2^\alpha = C_1$, then we can say that each
open ball $B(x,r)$ can be covered by a family of at most
$(2^l)^\alpha$ balls of radius $2^{-l} r$.  Note that if $(M, d(x,y))$
has the doubling property, then so does $(M, d(x,y)^a)$ for any
positive real number $a$.

	Suppose that $(M, d(x,y))$ is a quasimetric space, and that
$\mu$ is a nonnegative Borel measure on $M$.  Let us assume that open
balls in $M$ are Borel sets.  Of course, if $d(x,y)$ is a metric, then
open balls are open sets, and hence are Borel sets.  In practice, the
quasimetrics that one would consider do have this property, and anyway
one could make adjustments if necessary.  One says that $\mu$ is a
\emph{doubling measure} if the $\mu$-measure of open balls is positive
and finite, and if there is a positive real number $C_2$ such that
\begin{equation}
	\mu(B(x,2r)) \le C_2 \, \mu(B(x,r))
\end{equation}
for all $x \in M$ and $r > 0$.  A basic fact is that if there is a doubling
measure on $(M, d(x,y))$, then $(M, d(x,y))$ is doubling as a quasimetric
space.

	A quasimetric space $(M, d(x,y))$ equipped with a doubling
measure $\mu$ is often called a \emph{space of homogeneous type}.  As
in \cite{CW1, CW2}, a lot of real-variable methods in harmonic
analysis carry over to spaces of homogeneous type.  See \cite{eli1,
SW} for the classical setting of harmonic analysis on Euclidean
spaces, and see \cite{stevek, jean-lin, MS1, MS1, eli2} for more
information related to real-variable methods, doubling measures,
spaces of homogeneous type, etc.

\subsubsection*{Spaces at infinity of hyperbolic groups}

	Let $\Gamma$ be a finitely-presented group which is
hyperbolic.  Associated to $\Gamma$ is a space $\Sigma$ which is a
kind of ``space at infinity'' or ideal boundary of $\Gamma$,
consisting of equivalence classes of asymptotic directions in
$\Gamma$.  This space is a compact Hausdorff topological space of
finite dimension, as on p110-1 of \cite{misha2}, and it contains a
copy of the Cantor set as soon as it has at least three elements.  If
$\Sigma$ has at most two elements, then $\Gamma$ is said to be
\emph{elementary}.  For a free group with at least two generators the
space at infinity is homeomorphic to a Cantor set, while ${\bf Z}$, a
free group with one generator, is elementary and has two points in the
space at infinity.  If $\Gamma$ is the fundamental group of a closed
Riemann surface of genus at least $2$, then $\Sigma$ is homeomorphic
to the unit circle in ${\bf R}^2$.  More generally, if $\Gamma$ is the
fundamental group of a compact $n$-dimensional Riemannian manifold
without boundary with strictly negative curvature, then $\Sigma$ is
homeomorphic to the unit sphere ${\bf S}^{n-1}$ in ${\bf R}^n$.

	Actually, the space at infinity is defined for any hyperbolic
metric space in \cite{misha2}, and this can be specialized to a
hyperbolic group.  It is often preferable to work with metric spaces
which are ``geodesic'', in the sense that any pair of points can be
connected by a curve whose length is equal to the distance between the
two points.  It is often useful to think of a hyperbolic group as acting
on a geodesic hyperbolic metric space by isometries, and to use that
to study the space at infinity.

	It does not customarily seem to be said this way, but I think
it is fair to say that what are basically defined on the space at
infinity are quasimetrics, at least initially.  More precisely, it is
more like the logarithm of a quasimetric, or, in other words, there is
a one-parameter family of quasimetrics which are powers of each other.
In Section 7.2 of \cite{misha2} one takes a different route, in effect
compactifying a geodesic hyperbolic metric space by looking at
different measurements of lengths of curves which take densities into
account, densities which decay suitably at infinity.  For a parameter
in the density in an appropriate range, this measurement of lengths
of curves leads to a measurement of distance on the compactification
with nice properties, including upper and lower bounds by positive
constant multiples of quantities defined more directly.  By defining
distance in terms of lengths of curves in the compactification one
gets an actual metric in particular, i.e., with the usual triangle
inequality.

	It should perhaps be emphasized that in measuring distances
between points at infinity through weighted lengths of curves, the
curves are going through the hyperbolic metric space; curves in the
space at infinity are another matter, especially with about the
correct length.

	In nice situations, such as hyperbolic groups, and universal
coverings of compact Riemannian manifolds without boundary and
strictly negative curvature in particular, there are doubling
conditions on the space at infinity.  Compare with \cite{pierre}.
There are also interesting measures around, as in \cite{coornaert}.

	A well-known result of Borel \cite{borel, raghunathan} says
that simply-connected symmetric spaces can be realized as the
universal covering of a compact manifold.  If the symmetric space is
of noncompact type and rank $1$, it has negative curvature, and thus
the fundamental group of the compact quotient, which is a uniform
lattice in the group of isometries of the symmetric space, is a
hyperbolic group.  If the symmetric space is a classical hyperbolic
space of dimension $n$, with constant negative curvature, then the
space at infinity can be identified with a Euclidean sphere of
dimension $n - 1$.  If the symmetric space is a complex hyperbolic
space of complex dimension $m$, then the space at infinity can be
identified topologically with a Euclidean sphere of real dimension $2m
- 1$, but the geometry corresponds to a sub-Riemannian or
Carnot--Carath\'eodory space when $m \ge 2$, associated to a
distribution of hyperplanes in the tangent bundle of the sphere.  One
can think of the sphere as being the unit sphere in ${\bf C}^m$, and
the hyperplanes in the tangent bundle are the maximal complex
subspaces.  For other symmetric spaces of noncompact type and rank
$1$, one again obtains topological spheres of dimension $1$ less than
the real dimension of the symmetric space, and with sub-Riemannian
structures coming from distributions of planes of larger codimension
in the tangent bundle.

\end{document}